\documentclass[amsppt,11pt]{amsart}

\usepackage{amsmath}
\usepackage{amsthm}
\usepackage{amssymb}
\usepackage{verbatim}

\allowdisplaybreaks[4]

\def\R{\mathbb {R}}
\def\Z{\mathbb {Z}}

\def\P{\mathbb{P}}

\def\1{{\mathbf 1}}
\def\0{{\mathbf 0}}

\newtheorem{theorem}{Theorem}[section]

\def\supp{\mathrm{Supp\ }}

\title[Instability of random walks on groups]
{Instability of set recurrence and Green's function on groups with the
Liouville property}
\author[I. Benjamini\,\, D. Revelle]
{Itai Benjamini, David Revelle$^*$}
\date{October 1,2003.
\newline\indent
$^*$Research partially supported by an NSF postdoctoral
fellowship.}

\renewcommand{\k}{\mathbf {k}}
\renewcommand{\l}{\mathbf {l}}
\newcommand{\y}{\mathbf{y}}
\newcommand{\x}{\mathbf {x}}

\newcommand{\e}{\mathbf {e}}
\renewcommand{\u}{\mathbf {u}}
\renewcommand{\v}{\mathbf {v}}

\begin{document}

\begin{abstract}
Let $\mu$ and $\nu$ be probability measures on a group $\Gamma$ and let
$G_\mu$ and $G_\nu$ denote Green's function with respect to $\mu$ and
$\nu$.  The group $\Gamma$ is said to admit instability of Green's
function if there are symmetric, finitely supported measures $\mu$ and
$\nu$ and a sequence $\{x_n\}$ such that $G_\mu(e, x_n)/G_\nu(e,x_n)
\rightarrow 0$, and $\Gamma$ admits instability of recurrence if there is
a set $S$ that is recurrent with respect to $\nu$ but transient with
respect to $\mu$.
We give a number of examples of groups that have the Liouville property
but have both types of instabilities.  Previously known groups with these
instabilities did not have the Liouville property.
 \end{abstract}

\maketitle

\section{Introduction}

Let $\Gamma$ be an infinite graph.  Recall that simple
random walk on $\Gamma$ is a
Markov process such that $p(x,y)=1/deg(x)$ if $x \sim y$ and $0$
otherwise,
and a function $f$ is said to be harmonic if $f(x)=\sum p(x,y) f(y)$.
In the case when $\Gamma$ is the Cayley graph of a group, a common
generalization is to consider walks generated by a probability measure
$\mu$, i.e.,
$p(x,y)=\mu(x^{-1} y)$.  A number of properties of the long term behavior
of a random walk on a group are tied to the existence of non-constant bounded
harmonic functions.  A natural question to ask is to what
extent is the behavior of the random walk, or behavior of harmonic
functions,
affected by small perturbations in the graph $\Gamma$ or generating measure
$\mu$?

For metric spaces $\Gamma$ and $\Gamma^\prime$,
a map $\phi:
\Gamma \rightarrow \Gamma^\prime$
 is said to be a quasi-isometric
embedding
if there exist constants $\lambda \geq 1$ and $C\geq 0$ such that
$$
\frac{1}{\lambda} d(x,y) -C \leq d(\phi(x),\phi(y)) \leq \lambda d(x,y)
+C
$$
where $d(\cdot, \cdot)$ is the edge distance in the graph.  A
quasi-isometric embedding
$\phi$ is called a quasi-isometry if in addition
there is a constant $D\geq 0$ such that for any $y\in \Gamma^\prime$
there exists an $x\in \Gamma$ with $d(\phi(x),y)\leq D$.
An important example of a quasi-isometry comes from changing generating
sets for a group.
Let $S$ and $S^\prime$ be
two finite generating sets for a group and let $\Gamma$ and
$\Gamma^\prime$
be the metric
spaces whose elements are elements of the group and where the metric is
the world length with respect to $S$ and $S^\prime$ respectively.  These
spaces $\Gamma$ and $\Gamma^\prime$ are quasi-isometric.

Fix a basepoint $x_0 \in \Gamma$.  A set $S \subset \Gamma$
is said to be recurrent if simple random walk
starting at
$x_0$ visits $S$ infinitely often with probability 1.  When $\Gamma$ is a
group, we let $x_0$ be the identity, and say that $S$ is $\mu$-recurrent
if the walk generated by $\mu$ visits $S$ infinitely often with
probability 1.  Likewise, a set $S$ is said to be transient if simple
random walk almost surely visits $S$ only finitely often.

A function $f: \Gamma\rightarrow \R$ is $\mu$-harmonic if it satisfies the
mean value property, namely
$f(x)= \sum f(y)\mu(x^{-1}y)$.  On a general graph, $f$ is harmonic means
that $f(x)$ is the average of $f$ over the neighbors of $x$
The graph $\Gamma$
is said to have the Liouville property if all bounded harmonic functions
are constant.  Lyons \cite{TLyons87} took advantage of the fact
that the free group admits instability of recurrence to construct
a pair of quasi-isometric graphs such that one has the Liouville property
and the other one does not.

A second construction, due to
Benjamini and Schramm \cite{BenSch} is as follows:
suppose that $\Gamma$ is a Cayley graph and $S$ is $\mu$-recurrent but not
$\nu$-recurrent.  Let $\Phi:S\rightarrow \Z$ be any injection, and
consider the graph with vertex set $\Z^4 \cup \Gamma$ and edges consisting
of edges in $\Z^4$ and $\Gamma$, along with edges of the form $x \sim
(\Phi(x),0,0,0)$.  Because $S$ is $\mu$-recurrent and simple random walk
on $\Z^3$ is transient, a walk generated by $\mu$ will eventually be
absorbed by the $\Z^4$ component of the graph, and the fact that $\Z^4$ is
Liouville implies that the larger graph is as well.  A walk generated by
$\nu$ has a positive probability of not being absorbed by $\Z^4$, and
letting $f(y)$ be the probability of being absorbed in $\Z^4$ beginning at
$y$ gives a non-constant, bounded $\nu$-harmonic function.  More details
of the argument as well as a more general statement are in \cite{BenSch}.

Despite the existence, and in fact easy construction, of pairs of
quasi-isometric graphs such that one has the Liouville property while the
other does not, a major open question about random walks on groups is
whether or not having the Liouville property is an invariant for all
finitely supported, symmetric probability measures on a group.
There are some cases for which it is known that all such measures either
admit the Liouville property or do not: non-amenable groups are not
Liouville, while groups of polynomial volume growth are, as are polycyclic
groups.  The examples for which the Liouville property is not preserved
under quasi-isometry involve graphs which are not vertex transitive, and
it is possible that the vertex transitivity of Cayley graphs mean that the
Liouville property is preserved under changes of measure.  A recent survey
of some of these boundary properties is in \cite{KaiWoe}.

The nature of these examples may seem to suggest that groups
that admit an instability of recurrence should not have the
Liouville property.  This belief is perhaps supported by
the classical result that $\Z^d$ has the Liouville property and
 does not admit an instability of recurrence or Green's function.
Perhaps surprisingly, there are many other examples of groups which do.
The main focus of this paper is to present a number of examples
of groups with the Liouville property that admit instability
of both recurrence and Green's function.

For many natural examples, instability of recurrence and of Green's
function will go hand in hand.  The expected number of visits to a
recurrent set is, of course, infinite.  Conversely,
the arguments that we will use are
Borel-Cantelli based, and so the expected number of visits to the
transient sets we consider will be finite.

\begin{theorem}
\label{polycyclicthm}
Suppose that $G$ is a polycyclic group
of the form $G=\Z^m \ltimes \Z^d$, where for any $\k \in \Z^m$, the
eigenvalues of the matrix
associated with the action of $\k$ on $\Z^d$ are positive.
The group $G$ admits instability of both recurrence and Green's function
if and
only if $G$ has exponential volume growth.
\end{theorem}

At first, this theorem appears to be a very special case.  However,
a discrete subgroup $G$ of a Lie group $\Gamma$ contains a torsion free
finite index subgroup
$G^\prime$ that admits the exact sequence
$1 \rightarrow N \rightarrow G^\prime \rightarrow
\Z^d \rightarrow 1$, where $N$ is a torsion free
nilpotent group.
The main shortcomings of this theorem are thus the requirement that $N$
be abelian, and that the sequence splits as a semi-direct product.  These
restrictions appear as a result of the explicit computations
used in the proof and it is not clear that they should be necessary.

Another class of examples that admit instability of recurrence and Green's
function are wreath products over $\Z^d$.
Elements of the wreath product $G \wr H$
are pairs of the form
$(f,y)$, where $y\in H$ and $f\in \sum_{H} G$.  Recall that
elements of the direct sum $\sum_{H} G$ are functions from $H$
to $G$ which are equal to the identity for all but finitely many values.
The law is given by
$(f, y) (g, z)=(f s(y) g, yz)$, where $s(y)$
denotes translation by $y$, i.e. $s(y)g_i=g_{iy}$.
In this way, $H$ acts by
coordinate shift on $\sum_H G$ and $G \wr H=\sum_H G \rtimes H$ is the
resulting semi-direct product.

The case of random walk on $G=\Z_2$ and $H=\Z^d$ has received
particular attention because when $d=1$ or $2$, it is
an example of a group of exponential growth that has the Liouville
property,
and when $d\geq 3$, it is an example of an amenable group
which does not have the Liouville property.
These examples are often referred to as lamplighter groups because the
element $f\in \sum_{\Z^d} \Z_2$ can be viewed as a configuration of
lamps in $\Z^d$, finitely many of which are on, and $y\in \Z^d$ a
lamplighter who moves from lamp to lamp, turning them on and off.

\begin{theorem}
\label{wreaththm}
The lamplighter groups $G=H \wr \Z^d$ admit instability of recurrence and
Green's function.
\end{theorem}

Theorem \ref{wreaththm} is more surprising when $d=1$ or $2$ because in
those cases the group has the Liouville property.

\section{Proofs}

We begin by proving Theorem \ref{wreaththm} because
the proof contains the intuition behind the proof of
Theorem \ref{polycyclicthm} while
being technically simpler.

\begin{proof}[Proof of Theorem \ref{wreaththm}]
We first consider the case when $d=1$.
We need to produce two symmetric, finitely
supported probability
measures $\mu$ and $\nu$ and a set $S$ such that
$S$ is $\nu$-recurrent but not $\mu$-recurrent.
The proof will use the Borel-Cantelli lemmas and so will give
us the instability of Green's function for free.

The measure $\mu$ that we will consider is a convolution measure that
is typically used in random walks on lamplighter groups \cite{PitSC02}.
In the lamplighter description,
$\mu$ corresponds to the walk in which during
each time step the lamplighter randomizes the
current lamp, moves via simple random walk on $\Z$,
and then (for symmetry) randomizes the new lamp.
Formally, $\mu(f,x)=1/8$ if $x=\pm 1$ and
$f(i)=0$ for $i\notin \{0, x\}$, and $\mu(f,x)=0$ otherwise.

Let $\nu$ be the measure obtained by half
of the time sampling from $\mu$ and
half of the time having the lamplighter take
a step of size $8$ while not adjusting any
lamps.  Thus $\nu(f,x)=1/4$ if $x=\pm 8$ and
$f\equiv 0$, and
$\nu(f,x)=\mu(f,x)/2$ otherwise.

For an element $(f,x)\in \Z_2 \wr \Z$,
where $x\in \Z$ and $f\in \sum_{\Z} \Z_2$,
let $F_R$ be a flag
at the rightmost lamp that is on, or $0$ if no positive lamps are on,
i.e., $F_R=\max\{0, i: f(i)=1\}$.

For $\epsilon\in (0,2/5)$, let
  $S=S(\epsilon)$ be the set given by
$$
S=\left\{(f,x) : x-F_R \geq \frac{(2+5\epsilon) \log x}
{\log 2} , x\geq 1\right\}
$$

Let $(f_n,X_n)$ and $(g_n, Y_n)$
denote the walks generated by $\mu$ and $\nu$
respectively.  We will show that
with probability 1, $(f_n,X_n)\in S$ finitely often, and
$(g_n, Y_n)\in S$ infinitely often, i.e., $S$ is $\nu$-recurrent
but $\mu$-transient.

First we show that $(g_n,Y_n)\in S$ infinitely often.
Let $n_k=\inf\{n: Y_n \geq k^{1/2} \}$, $t_k=\lfloor \log k/\log
4\rfloor$,
and let $A_k$ be the event
$$
A_k=\{Y_{n_k+ t_k} = Y_{n_k}+ 8 t_k\}.
$$
In words, $A_k$ is the event that upon first reaching $k^{1/2}$, the
lamplighter then takes $t_k$ steps immediately to the right, reaching a
configuration in which the lamplighter is far enough to the right of
$F_R$ that the walk on the wreath product is in $S$.

 For large enough $k$, the events $A_k$ are independent because the time
intervals are disjoint.
Moreover,
whenever the event $A_k$ occurs, $x-F_R \geq 8 t_k$, and since $x$ is on
the order of $k^{1/2}$, $8 t_k$ is on the order of $4 \log_2 x$, which in
particular means that
 $(g_{n_k+t_k}, Y_{n_k+t_k})\in S$ once $k$ is large enough.
Since $\P(A_k)\geq 4^{-t_k}\geq k^{-1}$, infinitely many of the $A_k$
occur,
proving the first part of the claim.

The main point is to show that $(f_n,X_n)\in S$
finitely many times.
Let
$$M_n=\max_{k\leq n} X_k$$
be the rightmost lamp visited by the
lamplighter by time $n$.
Since $M_n$ is increasing, $M_n\geq X_n$ for all $n$, and $x-c\log x$ is
increasing for large $x$,
\begin{equation}
M_n-\frac{(2+5\epsilon) \log M_n}{\log 2} \geq X_n -\frac{(2+5\epsilon)
\log X_n}{\log 2}
\end{equation}
for all sufficiently large $n$.  In particular, it suffices to show that
the event
\begin{equation}
U_n=\left\{M_n-F_R(f_n,X_n) \geq \frac{(2+5\epsilon) \log M_n}{\log 2}
\right\}
\end{equation}
occurs finitely often.
The key is that our choice of generating set $\mu$
implies that each visited lamp is on with probability $1/2$, so
\begin{equation}
\label{maxineq}
\P[M_n-F_R\geq \alpha]=2^{-\alpha}.
\end{equation}
But
$M_n<n^{1/2-\epsilon}$ finitely often, so taking
$$\alpha=\lfloor (2+5\epsilon)\log_2 M_n \rfloor$$
yields a convergent series.  By the Borel-Cantelli lemma, finitely
many of the $U_n$ occur as required.

Moreover, this argument shows that the expected number of visits to
$S$ is finite.  Since the expected number of visits to $S$ is
$\sum_{x \in S} G_\mu(id, x)$ and the expected number of visits by a walk
generated by $\nu$ is infinite, this in turn implies that there is a
sequence $\{x_n\}\subset S$ with $G_\mu(id, x_n)/G_\nu(id, x_n)\rightarrow
0$, so Green's function is unstable as well.

Higher dimensions and the case when $H \neq \Z_2$ follow from a comparison
argument.  For the general case, take $\mu=\eta *\lambda * \eta$, where
$\eta$ is a symmetric measure satisfying $\eta(id)=1/2$ and $$\supp \eta
\subset \{(f,0): f(y)=e, \forall \ y\neq 0\},$$
that is to say $\eta$ performs a random walk in the current lamp with
holding probability $1/2$.  Take $\lambda$ to be a symmetric measure
with $\supp \lambda \subset \{(\e, \x)\}$ and $\sum_{\x\in \Z^d:
x_1=1}\lambda (\x) =1/2$, so $\lambda$ corresponds to the lamplighter
moving without adjusting any lamps in such a way that the projection onto
the first coordinate is simple random walk.  Take $\nu(f,\x)$ to be $1/4$
if $f=\e$ and $\x=(\pm 8, 0, \dots, 0)$ and
$\mu(f,\x)/2$ otherwise.  Let $F_R(f,x)$ be the rightmost point such that
$f(\y)\neq e$.  The behavior of $x-F_R$ for the walk generated by $\mu$ is
stochastically dominated by the behavior on $\Z_2 \wr \Z$, so $S$ is still
not $\mu$-recurrent, while the old argument still shows that $S$ is
$\nu$-recurrent.
 \end{proof}

\begin{proof}[Proof of Theorem \ref{polycyclicthm}]

We begin by assigning a coordinate system to $G=\Z^m\ltimes \Z^d$.
Because $G$ splits as a semi-direct product, we can uniquely write
any $g\in G$ as an ordered pair $(\k, \x)$, with $\k \in \Z^m$ and $\x\in
\Z^d$ such that the group multiplication law is given by
\begin{equation}
(\k,\x)(\l,\y)=(\k+\l, \x+\Psi_\k(\y))
\end{equation}
where by hypothesis
$\Psi_\k \in SL_d(\Z)$ has strictly positive eigenvalues.  Let
$\langle \e_1, \dots, \e_d \rangle$ be a basis for $\Z^d$, $\langle
\u_1, \dots, \u_m \rangle$ a basis for $\Z^m$, and let
$\lambda_{ij},
j=1, \dots, d$ be the eigenvalues for $\Psi_{\u_i}$.
Since $\Z^m$ is abelian, the matrices $\Psi_{\u_i}$ commute, so they
have a common eigenbasis $\langle \v_1, \dots, \v_d \rangle \in \R^d$.

Without loss of generality, we may assume that $\lambda_{11} \geq
\lambda_{ij}$ for all $i,j$.  The group $G$ has polynomial growth
whenever there is no dilation by the action.  Since our eigenvalues are
positive, this means that $G$ has polynomial growth iff
$\lambda_{ij}=1$ for all $i,j$.  In this case $G=\Z^{d+m}$, which is
classically known to have stability of recurrence and Green's function.
We will thus restrict our attention to the case when
$\lambda_{11}>1$, which is to say when $G$ has exponential volume growth.
We may also
assume that $(\e_i,\v_1)>0$ for all $i$, where $(\cdot,\cdot)$ denotes the
standard inner product.

For an element $g=(\k,\x)$, let
\begin{equation}
\phi(g)=\sum_{i=1}^m k_i \log \lambda_{i1}.
\end{equation}
The role of the functional $\phi$ will be analogous to the role of the
position of the lamplighter in the proof, and the analog to $F_R$ will be
played by $\log | (\x, \v_1)|$.
Denote the random walk $\xi_n$ by $\xi_n=(\k_n, \x_n)$, and let
$M_n=\max_{t\leq n} \phi(\xi_n)$.  Take
\begin{equation}
S=S(C)=\left\{ g: \phi(g)-\log |(\x,\v_1)| \geq C \log
\phi(g) \right\}.
\end{equation}

The two walks that we will consider on $G$ are given by the measures $\mu$
and $\nu$ as follows: let
$$\mu(\pm \u_i, \pm \e_j)=\frac{1}{6md}$$ and
$$\mu(\pm \u_i, 0)=\frac{1}{6m},$$
with $\mu(g)=0$ otherwise.
Define $\nu$ by letting $\nu(\pm \beta \u_i, 0)=\frac{1}{2} \frac{1}{2m}$
for a suitably large $\beta$ that we will choose later, and let
$\nu(g)=\mu(g)/2$ otherwise.

We will first choose $C$ large enough so that
$S$ is $\mu$-transient, and then choose $\beta$ depending on $C$
such that $S$ is $\nu$-recurrent.  Again, the Borel-Cantelli arguments
used to prove recurrence and transience will also prove the instability of
Green's function.

For the first part, we will prove that the analogy to (\ref{maxineq})
holds.  Let
$$t_0(n)=\min\{t\leq n: \phi(\k_t, \x_t)=M_n \}$$
and recursively define
$$t_\ell(n)=\max\{t\leq t_{\ell-1}(n): \phi(\xi_t) \leq
\phi(\xi_{t_{\ell-1}}) -\log 3 \}.$$

Since $M_n$ is the maximum value of a linear functional of the position
of a symmetric random walk on $\Z^m$, we have $M_n\leq
n^{1/2-\epsilon}$ finitely
often.
Let $y_\ell(n)=x_{t_\ell(n)+1}-x_{t_\ell(n)}$ be the increments in $x_n$
at time $t_\ell$.  By definition of $\mu$,
$$(y_\ell(n), \v_1)=\pm (\e_1,\v_1) \exp \phi(\xi_{t_\ell})$$
with probability $\frac{1}{3d}$ and $y_\ell=0$ with probability
$\frac{1}{3}$.  Moreover, this increment is independent of the change in
$\phi(\xi_{t})$.  This is analogous to the lamplighter randomizing the
current lamp independently of his movement in the previous proof.
 By the construction of the $t_\ell$,
$$\frac{\exp \phi(\xi_{t_{\ell+1}})}{\exp \phi(\xi_{t_{\ell}})} \geq 3
$$
and so
\begin{equation}
\label{maxineqsol}
\P[\log|(\x_n, \v_1)| < M_n-\beta(\log 3+ \lambda_{11})] \leq
\left(1-\frac{1}{3d}\right)^\beta.
\end{equation}
As before, this means that the events
\begin{equation}
U_n=\left\{M_n-\log|(\x_n, \v_1)| \geq C \log M_n \right\}
\end{equation}
occur finitely often for large enough $C$, and so for suitably large $C$,
the walk generated by $\mu$ is in $S(C)$ finitely often.

For the second part, let $n_r$ be the first time that $\phi(\xi_n)\geq
r-1$.  With probability at least $(4m)^{-\alpha}$,
$$\phi(\xi_{n_r + \alpha})-\phi(\xi_{n_r})\geq
\alpha \beta \log \lambda_{11}$$
while  $\x_{n_r}=\x_{n_r+\alpha}$.
But $(\x_{n_r},\v_1) \leq n_r \exp(r)$.  Thus, with probability at least
$(4m)^{-\alpha}$, we have the events
$$
\{ \log |(\x_{n_r+\alpha}, \v_1)| \leq \log n_r + r\}
$$
and
$$
\{\phi(\xi_{n_r+\alpha}) \geq r + \alpha \beta \log \lambda_{11} -1\}
$$
occurring simultaneously.  Taking $\alpha=\frac{\log r}{\log 4m}$ and
$\beta=\frac{6 C \log 4m}{\log n \log \lambda_{11}}$ does the job.
\end{proof}

\section{Further questions and remarks}

\begin{itemize}

\item
Instability of recurrence also allows us to construct examples in which
the rate of escape for a random walk varies under quasi-isometry.
For a random walk $\xi_n$ on a graph $\Gamma$, call $\liminf |\xi_n|/n$
the
speed of the random walk, where $|\xi_n|$ denotes the graph distance from
the starting point after $n$ steps.  When $\Gamma$ is the Cayley graph of
a group, $\lim |\xi_n|/n$ exists by subadditivity arguments.  In
the case
when $\mu$ is finitely supported and symmetric and $\xi_n$ is random walk
generated by $\mu$, then the speed is zero iff the group has the Liouville
property with respect to $\mu$.  For a survey of connections between
random walks and the Poisson boundary, see \cite{KaiWoe}.

 Suppose $\Gamma$ is the Cayley graph of a group for which the speed of
random walk is positive, such as $\Z_2 \wr \Z^3$.  Let $\mu$ and $\nu$ be
measures and $S$ a set such that $S$ is $\mu$-transient but
$\nu$-recurrent, and let the graph $G$ be obtained from $\Gamma$ by adding
an infinite ray at every point of $S$.  Because $S$ is $\mu$-transient,
the walk induced by $\mu$ on the new graph will still have positive speed,
but
because the expected return time on a ray is infinite, the walk induced by
$\nu$ will have zero speed.
\medskip

\item
Instability of recurrence is robust with respect to quotients.  More
precisely, note that if $N \lhd G$ and the quotient $G/N$ admits
instability of either recurrence or Green's function, then so does $G$.
To see this, suppose that $\mu$ and $\nu$ are measures on $G/N$ such that
$S$ is transient with respect to $mu$ and recurrent with respect to $\nu$.
Let $\tilde \mu$ and $\tilde \nu$ be measures on $G$ such that their
projections onto $G/N$ are $\mu$ and $\nu$.  Then $SN$ is transient with
respect to $\tilde \mu$ and recurrent with respect to $\tilde \nu$.

The next question is are these instabilities robust
with respect to extensions?  If $H \subset G$ and $H$ admits instability
of recurrence, does $G$ also admit instability of recurrence?
\medskip

\item
The main question is what is the key property of a group that
causes it to admit instability of recurrence or Green's function?  The
examples presented suggest that maybe the key is the volume growth of the
group.  Is it true
that a group admits instability of Green's function iff its volume growth
is exponential (or perhaps super-polynomial)?
Benjamini \cite{Ben91} gave examples of non-vertex transitive graphs with
polynomial growth that have instability of recurrence, so any equivalence
would need to use vertex transitivity.
Moreover, since the sets $S$ that we considered in the proofs of Theorems
\ref{polycyclicthm} and \ref{wreaththm} use the projection onto $\Z$ in
their descriptions, they
rely heavily on the fact that they contain elements of infinite order.  It
seems possible that
perhaps torsion groups of exponential growth have different behavior.

\end{itemize}

\bibliographystyle{amsplain}
\bibliography{../refs}

\providecommand{\bysame}{\leavevmode\hbox to3em{\hrulefill}\thinspace}
\providecommand{\MR}{\relax\ifhmode\unskip\space\fi MR }
\providecommand{\MRhref}[2]{%
  \href{http://www.ams.org/mathscinet-getitem?mr=#1}{#2}
}
\providecommand{\href}[2]{#2}
\begin{thebibliography}{1}

\bibitem{Ben91}
Itai Benjamini, \emph{Instability of the {L}iouville property for
  quasi-isometric graphs and manifolds of polynomial volume growth}, J.
  Theoret. Probab. \textbf{4} (1991), no.~3, 631--637. \MR{92f:31011}

\bibitem{BenSch}
Itai Benjamini and Oded Schramm, \emph{Harmonic functions on planar and almost
  planar graphs and manifolds, via circle packings}, Invent. Math. \textbf{126}
  (1996), no.~3, 565--587. \MR{97k:31009}

\bibitem{KaiWoe}
Vadim~A. Kaimanovich and Wolfgang Woess, \emph{Boundary and entropy of space
  homogeneous {M}arkov chains}, Ann. Probab. \textbf{30} (2002), no.~1,
  323--363. \MR{2003d:60152}

\bibitem{TLyons87}
Terry Lyons, \emph{Instability of the {L}iouville property for quasi-isometric
  {R}iemannian manifolds and reversible {M}arkov chains}, J. Differential Geom.
  \textbf{26} (1987), no.~1, 33--66. \MR{88k:31012}

\bibitem{PitSC02}
C.~Pittet and L.~Saloff-Coste, \emph{On random walks on wreath products}, Ann.
  Probab. \textbf{30} (2002), no.~2, 948--977. \MR{2003d:60013}

\end{thebibliography}

\end{document}